\documentclass[10pt,leqno]{article}
\usepackage{amssymb,amsbsy,amsmath,amsfonts,amssymb,amscd, mathrsfs}
\usepackage{graphicx} 

\usepackage[all]{xy} 
\usepackage{color}

\usepackage[english]{babel}
\usepackage[T1]{fontenc}
\usepackage{indentfirst}

\makeatletter
\@addtoreset{equation}{section}
\makeatother

\newtheorem{statement}{}[section]
\newtheorem{theorem}[statement]{Theorem}
\newtheorem{lemma}[statement]{Lemma} 

\newtheorem{proposition}[statement]{Proposition}

\newcommand\C{\mathbb C}

\newcommand\T{\mathbb T}
\newcommand\D{\mathbb D}

\newcommand\e{{\rm e}}

\newcommand\eps{\varepsilon}
\newcommand\ind{{\rm 1\kern-.30em I}}
\newcommand\qed{\hfill $\square$}
\renewcommand \Re{{\mathfrak R}{\rm e}\,}
\renewcommand \Im{{\mathfrak I}{\rm m}\,}

\let\phi=\varphi

\title{\bf Two remarks on  composition operators on the  Dirichlet space}
\author{\it Daniel Li, Herv\'e Queff\'elec, \\ 
\it Luis Rodr{\'\i}guez-Piazza\footnote{Supported by a Spanish research project MTM 2012-05622.}}

\date{\footnotesize \today}

\begin{document}

\maketitle 

\noindent{\bf Abstract.} \emph{We show that the decay of approximation numbers of compact composition operators on the Dirichlet space $\mathcal{D}$ can 
be as slow as we wish. We also prove the optimality of a result of O.~El-Fallah, K.~Kellay, M.~Shabankhah and H.~Youssfi 
on boundedness on $\mathcal{D}$ of self-maps of the disk all of whose powers are norm-bounded in $\mathcal{D}$.}
\medskip

\par\medskip



\section{Introduction}

Recall that if $\varphi$ is an analytic self-map of $\D$, a so-called {\it Schur function}, the composition operator $C_\varphi$ associated to $\varphi$ is formally 
defined by
\begin{displaymath} 
C_{\varphi}(f)=f\circ \varphi \, .
\end{displaymath} 
The  Littlewood subordination principle (\cite{COWMAC}, p.~30) tells us that $C_\varphi$ maps the Hardy space $H^2$ to itself for every Schur function 
$\varphi$. Also recall that if $H$ is a Hilbert space and  $T \colon H \to H$ a bounded linear operator, the $n$-th approximation number $a_{n}(T)$ of $T$ is 
defined as 
\begin{equation} \label{nbres approx}
\qquad a_{n}(T) = \inf\{\Vert T - R \Vert \, ; \ \text{rank}\, R  < n \},  \quad  n = 1, 2, \ldots \, .
\end{equation} 
In \cite{JAT}, working on that  Hardy space $H^2$ (and also on some weighted Bergman spaces), we have undertaken the study of approximation numbers 
$a_{n} (C_\varphi)$ of composition operators $C_\varphi$, and proved among other facts the following:
\begin{theorem} \label{JAT} 
Let $(\varepsilon_n)_{n \geq 1}$ be a non-increasing sequence of positive numbers tending to $0$. Then, there exists a compact composition operator 
$C_\varphi$ on $H^2$ such that 
\begin{displaymath} 
\liminf_{n \to \infty} \frac{a_{n} (C_\varphi)}{\varepsilon_n} > 0 \, .
\end{displaymath} 
As a consequence, there are composition operators on $H^2$ which are compact but in no Schatten class.
\end{theorem}

The last item had been previously proved by Carroll and Cowen (\cite{CARCOW}), the above statement with approximation numbers being more precise. \par
\smallskip

For the Dirichlet space, the situation is more delicate because not every analytic self-map of $\D$ generates a bounded composition operator on $\mathcal{D}$.
When this is the case,  we will say that $\varphi$ is a \emph{symbol} (understanding ``of ${\mathcal D}$''). Note that every symbol is necessarily in ${\cal D}$. 
\par

In \cite{PDHL}, we have performed a similar study on that Dirichlet space $\mathcal{D}$, and established  several results on approximation numbers in that new 
setting, in particular the existence of symbols $\varphi$ for which $C_\varphi$ is compact without being in any Schatten class $S_p$. But we have not been able in 
\cite{PDHL} to prove a full analogue of Theorem~\ref{JAT}. Using a new approach, essentially based on Carleson embeddings and the Schur test, we are now 
able to prove that analogue.

\begin{theorem}\label{NEW} 
For every sequence $(\varepsilon_n)_{n \geq 1}$ of positive numbers tending to $0$, there exists a compact composition operator 
$C_\varphi$ on the Dirichlet space $\mathcal{D}$  such that 
\begin{displaymath} 
\liminf_{n\to \infty} \frac{a_{n} (C_\varphi)}{\varepsilon_n} > 0 \, .
\end{displaymath} 
\end{theorem}
\medskip

Turning now to the question of necessary or sufficient conditions for a Schur function  $\varphi$ to be a symbol, we can observe that, since 
$(z^n / \sqrt n)_{n\geq 1}$ is an orthonormal sequence in $\mathcal{D}$ and since formally $C_{\varphi} (z^n) = \varphi^n$, a necessary condition  is as 
follows:
\begin{equation}\label{nec} 
\varphi \text{ is a symbol} \quad \Longrightarrow \quad \Vert \varphi^n \Vert_{\mathcal{D}} = O\, (\sqrt n) \, .
\end{equation}
It is worth noting that, for any Schur function, one has:
\begin{displaymath} 
\varphi \in \mathcal{D} \quad \Longrightarrow  \quad \Vert \varphi^n \Vert_{\mathcal{D}} = O\, (n)
\end{displaymath} 
(of course, this is an equivalence). Indeed, anticipating on the next section, we have for any integer $n \geq 1$:
\begin{align*}
\Vert \varphi^n \Vert_{\mathcal{D}}^2 
& = |\varphi (0)|^{2 n}+\int_{\D} n^2 \,|\varphi (z)|^{2 (n - 1)}| \varphi ' (z)|^2 \, dA (z) \\
& \leq |\varphi (0)|^{2} + \int_{\D} n^2 \, |\varphi ' (z)|^2 \, dA (z) 
\leq n^2 \Vert \varphi \Vert_{\mathcal{D}}^{2}, 
\end{align*}
giving the result.\par 

Now, the following sufficient condition was given in \cite{EKSY}:
\begin{equation}\label{suf} 
\Vert \varphi^n\Vert_{\mathcal{D}} = O\, (1) \quad \Longrightarrow \quad \varphi \text{ is a symbol} \, .
\end{equation}
In view of \eqref{nec}, one might think of improving this condition, but it turns out to be optimal, as says the second main result of that paper.

\begin{theorem}\label{OPT} 
Let $(M_n)_{n \geq 1}$ be an arbitrary sequence of positive numbers tending to $\infty$. Then, there exists a Schur function $\varphi \in \mathcal{D}$ such that:
\par\smallskip 

1) $\Vert \varphi^{n} \Vert_{\mathcal{D}} = O\, (M_n)$ as $n\to \infty$; \par
\smallskip

2) $\varphi$ is not a symbol on $\mathcal{D}$.
\end{theorem}

The organization of that paper will be as follows: in Section~\ref{notations}, we give the notation and background. In Section~\ref{proof Th 1}, we prove 
Theorem~\ref{NEW}; in Section~\ref{proof Th 2}, we  prove Theorem~\ref{OPT}; and we end with a section of remarks and questions.


\section{Notation and background.} \label{notations}

We denote by $\D$ the open unit disk of the complex plane and by $A$ the normalized area measure $dx \, dy / \pi$ of $\D$. The unit circle is denoted by 
$\T = \partial \D$. The notation $A \lesssim B$ indicates that $A \leq c \, B$ for some positive constant $c$.
\par \smallskip

A Schur function is an analytic self-map of $\D$ and the associated composition operator is defined, formally, by $C_\phi (f) = f \circ \phi$. The operator 
$C_\varphi$ maps the space ${\cal H}{ol}\, (\D)$  of holomorphic functions on $\D$ into itself.
\par\smallskip

The Dirichlet space $\mathcal{D}$ is the space of analytic functions $f \colon \D \to \C$ such that 
\begin{equation} 
\| f \|_{\cal D}^2 := | f (0)|^2 + \int_\D |f ' (z) |^2 \, dA (z) < + \infty \, .
\end{equation} 
If $f (z) = \sum_{n = 0}^\infty c_n z^n$, one has:
\begin{equation} 
\| f \|_{\cal D}^2 = |c_0|^2 + \sum_{n = 1}^\infty n \, |c_n|^2 \, .
\end{equation} 
Then $\| \, . \, \|_{\cal D}$ is a norm on ${\cal D}$, making ${\cal D}$ a Hilbert space, and $\| \, . \, \|_{H^2} \leq \| \, . \, \|_{\cal D}$. For further 
information on the Dirichlet space, the reader may see \cite{survey} or \cite{Ross}. \par \smallskip

The Bergman space ${\mathfrak B}$ is the space of analytic functions $f \colon \D \to \C$ such that: 
\begin{displaymath} 
\| f \|_{\mathfrak B}^2 := \int_\D |f (z) |^2 \, dA (z) < + \infty \, .
\end{displaymath} 
If $f (z) = \sum_{n = 0}^\infty c_n z^n$, one has $\| f \|_{\mathfrak B}^2 = \sum_{n = 0}^\infty \frac{|c_n|^2}{n + 1}$. If $f \in \mathcal{D}$, one has 
by definition:
\begin{displaymath}
\| f \|_{\mathcal{D}}^2 = \| f ' \|_{\mathfrak B}^2 +|f (0)|^2 \,.
\end{displaymath}

Recall that, whereas every Schur function $\phi$ generates a bounded composition operator $C_\phi$ on Hardy and Bergman spaces, it is no longer the case for the 
Dirichlet space (see \cite{McCluer-Shapiro}, Proposition~3.12, for instance). \par
\medskip 

We denote by $b_n (T)$ the $n$-th \emph{Bernstein number} of the operator $T \colon H\to H$, namely:  
\begin{equation} \label{Bernstein} 
b_n (T) = \sup_{\dim E = n} \Big( \inf_{f \in S_E} \Vert T x \Vert\Big)  \, 
\end{equation}
where $S_E$ denotes the unit sphere of $E$. It is easy to see (\cite{PDHL}) that 
\begin{displaymath} 
\qquad b_{n} (T) = a_{n} (T) \quad \text{for all } n \geq 1 \, .
\end{displaymath} 
(recall that the approximation numbers are defined in \eqref{nbres approx}). \par\smallskip

 If $\varphi$ is a Schur function, let 
\begin{equation} 
n_{\varphi} (w) = \# \{z\in \D \, ; \ \varphi (z) = w\} \geq 0 
\end{equation} 
be the associated \emph{counting function}. If $f \in \mathcal{D}$ and $g = f\circ \varphi$, the change of variable formula provides us with the useful following 
equation (\cite{ZOR}, \cite{PDHL}):
\begin{equation}\label{utile} 
\quad \int_{\D} |g'(z)|^2 \, dA (z) = \int_{\D} |f '(w)|^2 \,n_{\varphi} (w) \, dA (w)  
\end{equation} 
(the integrals might be infinite). In those terms, a necessary and sufficient condition for $\varphi$ to be a symbol is as follows (\cite{ZOR}, Theorem~1). Let:
\begin{equation} \label{def-ninaz} 
\rho_{\varphi}(h) = \sup_{\xi \in \T} \int_{S (\xi, h)} n_{\varphi} \, dA
\end{equation} 
where $S (\xi, h) = \D \cap D (\xi, h)$ is the Carleson window centered at $\xi$ and of size $h$. Then $\phi$ is a symbol if and only if:
\begin{equation}\label{ninaz}
\sup_{0 < h < 1}\frac{1}{h^2}\,\rho_{\varphi} (h) <\infty.
\end{equation}
This is not difficult to prove. In view of \eqref{utile}, the boundedness of $C_\varphi$ amounts to the existence of a constant $C$ such that:
\begin{displaymath} 
\int_{\D} |f ' (w)|^2 \, n_{\varphi} (w) \, dA (w) \leq C \int_{\D} |f '(z)|^2 \, dA (z) \, , \quad \forall f \in \mathcal{D}.
\end{displaymath} 
Since $f ' = h$ runs over $\mathfrak{B}$ as $f$ runs over $\mathcal{D}$, and with equal norms, the above condition reads:
\begin{displaymath} 
\int_{\D} |h (w)|^2 \, n_{\varphi} (w) \, dA (w) \leq C \int_{\D} |h (z)|^2 \, dA (z) \, , \quad \forall h\in \mathfrak{B}.
\end{displaymath} 
This exactly means that the measure $n_\varphi\, dA$ is a Carleson measure for $\mathfrak{B}$. Such measures have been characterized in \cite{HAS} and that 
characterization gives \eqref{ninaz}. \par\smallskip
 
But this condition is very abstract and difficult to test, and sometimes more ``concrete'' sufficient conditions are desirable. In \cite{PDHL}, we proved that, even 
if the Schur function extends continuously to $\overline{\D}$, no Lipschitz condition of order $\alpha$, $0 < \alpha < 1$,  on $\varphi$ is sufficient for ensuring 
that $\varphi$ is a symbol. It is worth noting that the limiting case $\alpha = 1$, so restrictive it is, guarantees the result.

\begin{proposition} 
Suppose that the Schur function $\varphi$ is in the analytic Lipschitz class on the unit disk, i.e. satisfies:
\begin{displaymath} 
\qquad | \phi (z) - \phi (w) | \leq C \, |z - w| \, , \quad \forall z, w \in \D\, .
\end{displaymath} 
Then $C_\varphi$ is bounded on $\mathcal{D}$.
\end{proposition}
\noindent {\bf Proof.} Let $f \in \mathcal{D}$; one has:
\begin{align*}
\Vert C_{\varphi} (f) \Vert_{\mathcal{D}}^2 
& =| f \big( \phi (0) \big)|^2 + \int_{\D} \vert f ' \big(\varphi (z) \big) \vert^{2} \vert \varphi ' (z) \vert^2 \, dA (z) \\ 
& \leq | f \big( \phi (0) \big)|^2 + \Vert \varphi ' \Vert_\infty^2 \int_{\D} \vert f ' \big(\varphi (z) \big) \vert^{2} \, dA (z) \, . 
\end{align*}
This integral is nothing but $\|C_\phi (f ') \|_{\mathfrak B}^2$ and hence, since $C_\varphi$ is bounded on the Bergman space $\mathfrak{B}$, we have, 
for some constant $K_1$:
\begin{displaymath} 
\int_{\D} \vert f ' \big(\varphi (z) \big) \vert^{2} \, dA (z) \leq K_1^2 \| f ' \|_{\mathfrak B}^2  \leq K_1^2 \| f \|_{\cal D}^2 \, .
\end{displaymath} 
On the other hand, 
\begin{displaymath} 
|f \big( \phi (0) \big)| \leq ( 1 - |\phi (0)|^2)^{- 1 / 2} \| f \|_{H^2} \leq ( 1 - |\phi (0)|^2)^{- 1 / 2} \| f \|_{\cal D} \, , 
\end{displaymath} 
and we get 
\begin{displaymath} 
\Vert C_{\varphi} (f) \Vert_{\mathcal{D}}^2 \leq K^2  \| f \|_{\cal D}^2 \, ,
\end{displaymath} 
with $K^2 = K_1^2 + ( 1 - |\phi (0)|^2)^{- 1}$. \qed 

\goodbreak

\section{Proof of Theorem \ref{NEW}} \label{proof Th 1}

We are going to prove Theorem~\ref{NEW} mentioned in the Introduction, which we recall here. 
\begin{theorem}\label{slow} 
For every sequence $(\varepsilon_n)$ of positive numbers with limit $0$, there exists a compact composition operator $C_{\varphi}$ on $\mathcal{D}$ such that 
\begin{displaymath} 
\liminf_{n\to \infty} \frac{a_{n} (C_\varphi)}{\varepsilon_n} > 0 \, .
\end{displaymath} 
\end{theorem}
\smallskip

Before entering really in the proof, we may remark that, without loss of generality, by replacing $\eps_n$ with $\inf (2^{- 8}, \sup_{k \geq n} \eps_k)$, we can, 
and do, assume that $(\varepsilon_n)_n$ decreases and $\varepsilon_1 \leq 2^{- 8}$. \par 

Moreover, we can assume that $(\eps_n)_n$ decreases ``slowly'', as said in the following lemma.

\begin{lemma} \label{tilde} 
Let $(\varepsilon_i)$ be a decreasing sequence with limit zero and let $0 < \rho < 1$. Then, there exists another sequence $(\widehat{\varepsilon_i})$, decreasing  
with limit zero, such that $\widehat{\varepsilon_i} \geq \varepsilon_i$ and  $\widehat{\varepsilon_{i + 1}} \geq \rho \, \widehat{{\varepsilon_i}}$, for every 
$i \geq 1$.
\end{lemma}

\noindent {\bf Proof.} We define inductively  $\widehat{\varepsilon_i}$ by  $\widehat{\varepsilon_1} = \varepsilon_1$ and 
\begin{displaymath} 
\widehat{\varepsilon_{i + 1}} = \max (\rho \, \widehat{\varepsilon_i}, \varepsilon_{i+1}).
\end{displaymath} 
It is seen by induction that  $\widehat{\varepsilon_i} \geq \varepsilon_i$ and that $\widehat{\varepsilon_i}$ decreases to a limit $a \geq 0$. If 
$ \widehat{\varepsilon_i} = \varepsilon_i$ for infinitely many indices $i$, we have $a = 0$. In the opposite case, 
$\widehat{\varepsilon_{i + 1}} = \rho \, \widehat{\varepsilon_i}$ from some index $i_0$ onwards, and again $a = 0$ since $\rho < 1$. \qed
\medskip

We will take $\rho = 1/2$ and assume for the sequel that $\eps_{i + 1} \geq \eps_i / 2$. \par \medskip

\noindent{\bf Proof of Theorem~\ref{slow}.} We first construct a subdomain $\Omega = \Omega_\theta$ of $\D$ defined by a cuspidal inequality:
\begin{equation}\label{pidal} 
\Omega = \{z = x + i y \in\D \, ; \ \vert y \vert < \theta (1 - x) \, , \  0 < x < 1 \} \, ,
\end{equation}
where $\theta \colon [0, 1] \to [0, 1[$ is a continuous increasing function such that 
\begin{equation}\label{theta}
\theta (0) = 0 \quad \text{and} \quad  \theta (1 - x) \leq 1 - x \, .
\end{equation}

Note that since $1 - x \leq \sqrt{1 - x^2}$, the condition $|y| < \theta (1 - x)$ implies that $z = x + i y \in \D$. Note also that $1 \in \overline{\Omega}$ and 
that $\Omega$ is a Jordan domain. \par\smallskip

We introduce a parameter $\delta$ with $\eps_1 \leq \delta \leq 1 - \eps_1$. We put:
\begin{equation}\label{put} 
\theta (\delta^j) = \eps_j \, \delta^j 
\end{equation} 
and  we extend  $\theta$ to an increasing continuous function from $(0, 1)$ into itself (piecewise linearly, or more smoothly, as one wishes). We  claim that:
\begin{equation} \label{nina} 
\theta (h) \leq h \quad \text{and} \quad  \theta (h) = o\, (h) \hbox{ as } h \to 0 \, .
\end{equation}
Indeed,  if $\delta^{j + 1} \leq h < \delta^j$, we have $\theta (h) / h \leq \theta (\delta^j) / \delta^{j + 1} = \eps_j / \delta$, which is 
$\leq \eps_1 / \delta \leq 1$ and which tends to $0$ with $h$. \par \smallskip

We define now $\phi = \phi_\theta \colon \overline{\D} \to \overline{\Omega}$ as a continuous map which is a Riemann map from $\D$ onto $\Omega$, 
and with $\varphi (1) = 1$ (a cusp-type map). Since $\phi$ is univalent, one has $n_\phi = \ind_\Omega$, and since $\Omega$ is bounded, $\phi$ defines a 
symbol on $\mathcal{D}$, by \eqref{ninaz}. Moreover, \eqref{nina} implies that $A [S (\xi, h) \cap \Omega] \leq h \, \theta (h)$ for every $\xi \in \T$; hence, 
$\rho_\varphi$ being defined in \eqref{def-ninaz}, one has $\rho_{\varphi} (h) = o \, (h^2)$ as $h \to 0^{+}$. In  view of \cite{ZOR}, this little-oh condition 
guarantees the compactness of $C_{\varphi} \colon \mathcal{D} \to \mathcal{D}$. \par
\smallskip

It remains to minorate its approximation numbers.\par 

The measure $\mu = n_{\varphi} \, dA$ is a Carleson measure for the Bergman space ${\mathfrak B}$, and it was proved in \cite{Dirichlet} that 
$C_{\varphi}^\ast C_\varphi$ is unitarily equivalent to the Toeplitz operator $T_\mu = I_{\mu}^\ast I_\mu  \colon {\mathfrak B} \to {\mathfrak B}$ 
defined by:
\begin{equation} \label{toep} 
\qquad T_{\mu} f (z) = \int_{\D} \frac{f (w)}{(1 - \overline{w} z)^2} \, dA (w) = \int_{\D} f (w) K_{w} (z) \, dA (w) \, , 
\end{equation}
where $I_\mu \colon {\mathfrak B} \to L^2 (\mu)$ is the canonical inclusion and $K_w$ the reproducing kernel of $\mathfrak{B}$ at $w$, i.e.  
$K_{w} (z )= \frac{1}{(1 - \overline{w} z)^2}\,$. \par

Actually, we can get rid of the analyticity constraint in considering,  instead of $T_\mu$, the operator 
$S_\mu = I_\mu I_{\mu}^\ast \colon L^{2} (\mu) \to L^{2} (\mu)$, which corresponds to the arrows:
\begin{displaymath} 
L^{2} (\mu) \mathop{\longrightarrow}^{I_\mu^\ast} \mathfrak{B} \mathop{\longrightarrow}^{I_\mu} L^{2} (\mu) \, .
\end{displaymath} 
We use the relation \eqref{toep} which implies:
\begin{equation} \label{equal}  
a_{n} (C_\varphi) = a_{n} (I_\mu) = a_{n} (I_{\mu}^\ast) = \sqrt{a_{n} (S_\mu)} \, .
\end{equation}

We set:
\begin{equation} \label{r_j} 
c_j = 1 - 2 \delta^j \quad \text{and} \quad r_j = \eps_j \, \delta^j 
\end{equation} 
One has $r_j = \eps_j (1 - c_j)/ 2$. 
\begin{lemma} \label{inclus} 
The disks $\Delta_j = D (c_j, r_j)$, $j \geq 1$, are disjoint and contained in $\Omega$.
\end{lemma}
\noindent{\bf Proof.} If $z = x + i y \in \Delta_j$, then $1 - x > 1 - c_j - r_j = (1 - c_j) (1 - \eps_j/2) = 2 \delta^j  (1 - \eps_j/2) \geq \delta^j$ and 
$|y| < r_j = \theta ( \delta^j)$; hence $\vert y \vert < \theta (\delta^j) \leq \theta (1 - x)$ and $z \in \Omega$. On the other hand, 
$c_{j + 1} - c_j = 2 (\delta^j - \delta^{j + 1}) = 2 (1 - \delta) \delta^j \geq 2 \eps_1 \delta^j \geq 2 \eps_j \delta^j = 2 r_j > r_j + r_{j + 1}$; hence 
$\Delta_j \cap \Delta_{j + 1} = \emptyset$. \qed 
\par \medskip

We will next  need a description of $S_\mu$.

\begin{lemma}\label{adj} 
For every $g \in L^2 (\mu)$ and every $z \in \D$:
\begin{align}\label{oint} 
I_{\mu}^\ast g (z) & = \int_\Omega \frac{g (w)}{(1 - \overline{w} z)^2} \, dA (w) \\
S_{\mu} g (z) & = \bigg(\int_\Omega \frac{g (w)}{(1 - \overline{w} z)^2} \, dA (w) \bigg)\, \ind_\Omega (z) \, .
\end{align}
\end{lemma}

\noindent{\bf Proof.} $K_w$ being the reproducing kernel of ${\mathfrak B}$, we have for any pair of functions $f \in \mathfrak{B}$ and $g \in L^2 (\mu)$:
\begin{align*}
\langle I_{\mu}^{*} g, f \rangle_{{\mathfrak B}} =\langle g, I_{\mu} f \rangle_{L^{2}(\mu)} 
& = \int_{\Omega} g (w) \overline{f (w)} \, dA (w) =\int_{\Omega} g (w) \, \langle K_w, f \rangle_{{\mathfrak B}} \, dA (w) \\ 
& =\big\langle \int_{\Omega} g (w) K_w \, dA (w), f \big\rangle_{{\mathfrak B}}, 
\end{align*}
so that  $I_{\mu}^{*} g = \int_{\Omega} g (w) K_w \, dA (w)$, giving the result. \qed
\medskip

In the rest of the proof, we fix a positive integer $n$ and put:
\begin{equation}\label{adopt} 
\qquad\qquad\qquad\quad  f_j = \frac{1}{r_j} \, \ind_{\Delta_j} \, , \qquad\quad  j = 1, \ldots, n\, . 
\end{equation}
Let:
\begin{displaymath} 
E = {\rm span}\, ( f_1, \ldots, f_n) \, .
\end{displaymath} 
This is an $n$-dimensional subspace of $L^2 (\mu)$. \par\smallskip

 The $\Delta_j$'s being disjoint, the sequence $(f_1, \ldots, f_n)$ is orthonormal in $L^{2}(\mu)$. Indeed, those functions have disjoint supports, so are 
orthogonal, and:
\begin{displaymath} 
\int\, f_j^2 \, d\mu = \int f_j^2 \, n_{\varphi}\, dA = \int_{\Delta_j} \frac{1}{r_j^2} \, dA = 1 \, .
\end{displaymath} 

We now estimate from below the Bernstein numbers of $I_{\mu}^{*}$. To that effect, we compute the scalar products 
$m_{i, j} = \langle I_{\mu}^{*} (f_i), I_{\mu}^{*} (f_j) \rangle$. One has:
\begin{align*}
m_{i, j} 
& =\langle f_i, S_{\mu} (f_j) \rangle = \int_{\Omega} f_i (z) \overline{S_{\mu} f_j (z)} \, dA (z) \\ 
& =\iint_{\Omega \times\Omega} \frac{f_i (z) \overline{f_j (w)}}{(1 - w \overline{z})^2} \, dA (z) \, dA (w) \\ 
& =\frac{1}{r_i r_j} \iint_{\Delta_i \times \Delta_j} \frac{1}{(1 - w \overline{z})^2} \, dA (z) \, dA (w) \, .
\end{align*}
\par\smallskip

\begin{lemma} \label{tec} 
We have 
\begin{equation} \label{dimanche} 
\qquad  m_{i, i} \geq \frac{\varepsilon_{i}^2}{32}, \qquad \text{and} \quad  
| m_{i, j}| \leq \varepsilon_i \, \varepsilon_j \,\delta^{j - i} \quad \text{for } i < j \, . 
\end{equation}
\end{lemma}

\noindent {\bf Proof.}  Set $\eps'_i = \frac{r_i}{1 - c_i^2} = \frac{\eps_i}{2 (1 + c_i)}$. One has $\frac{\eps_i}{4} \leq \eps'_i \leq \frac{\eps_i}{2}$. 
We observe that (recall that $A (\Delta_i) = r_i^2$):
\begin{displaymath} 
m_{i, i} - {\eps'_i}^{2}
= \frac{1}{r_i^2} \iint_{\Delta_i \times \Delta_i} \bigg[ \frac{1}{(1 - w \overline{z})^2} - \frac{1}{(1 - c_{i}^2)^2} \bigg] \, dA (z) \, dA (w) \, . 
\end{displaymath} 
Therefore, using  the fact that, for $z\in \Delta_i$ and $w \in \mathbb{D}$:
\begin{displaymath} 
| 1 - w \overline{z}| \geq 1 - | z | \geq 1 - c_i - r_i = 1 - c_i - \varepsilon_i \Big(\frac{1 - c_i}{2}\Big) \geq (1 - c_i) \Big(1 - \frac{\varepsilon_i}{2} \Big) 
\geq \frac{1 - c_i}{2}
\end{displaymath} 
and then the mean-value theorem, we get:
\begin{align*}
| m_{i, i} - {\eps'_i}^2 |  
& \leq   \frac{1}{r_i^2} \iint_{\Delta_i \times \Delta_i} \bigg| \frac{1}{(1 - w \overline{z})^2} - \frac{1}{(1 - c_i^2)^2} \bigg| \, dA (z) \, dA (w) \\ 
& \leq \frac{1}{r_i^2} \iint_{\Delta_i \times \Delta_i} \frac{32 \, r_i}{(1 - c_i)^3} \, dA (z) \, dA (w) \\ 
& = \frac{32 \, r_i^3}{(1 - c_i)^3} \leq 32 \times 8 \, {\varepsilon'_i}^3 \leq \frac{{\eps'_i}^2}{2} \, \raise 1pt \hbox{,} 
\end{align*}
since $\eps_i \leq \eps_1 \leq 2^{ - 8}$ implies that $\eps'_i \leq 1/ (32 \times 16)$. This gives us the lower bound 
$m_{i, i} \geq {\eps'_i}^2 /2 \geq \eps_i^2 / 32$.  \par 

Next, for $i < j$: 
\begin{align*}
| m_{i, j} | 
& \leq \frac{1}{r_i r_j} \iint_{\Delta_i \times \Delta_j} \bigg| \frac{1}{(1 - w \overline{z})^2} \bigg| \,  dA (z) \, dA (w) 
\leq \frac{1}{r_i r_j} \frac{4}{(1 - c_i)^2} r_i^2 r_j^2 \\ 
& = \frac{4 \, \varepsilon_i \, \varepsilon_j \, \delta^{i + j}}{4 \, \delta^{2 i}} = \varepsilon_i \, \varepsilon_j \, \delta^{j - i} \, ,
\end{align*}
and that ends the proof of  Lemma~\ref{tec}. \qed
\bigskip

We further write the $n \times n$ matrix $M = (m_{i, j})_{1 \leq i, j \leq n}$ as $M = D + R$ where $D$ is the diagonal matrix $m_i = m_{i, i}$ with 
$m_i \geq \frac{\eps_i^2}{32}$, $1 \leq i \leq n$. Observe that $M$ is nothing but the matrix of $S_\mu$ on the orthonormal basis $(f_1, \ldots, f_n)$ of $E$, 
so that we can identify $M$ and $S_\mu$ on $E$. \par

Now the following lemma will end the proof of Theorem~\ref{slow}. 
  
\begin{lemma}\label{cle} 
If $\delta \leq 1/200$, we have:
\begin{equation} \label{schur} 
\| D^{- 1} R \| \leq 1/2 \, .
\end{equation}
\end{lemma}

Indeed, by the ideal property of Bernstein numbers, Neumann's lemma and the relations:
\begin{displaymath} 
M = D (I + D^{- 1} R) \, , \quad  \text{and} \quad D = M Q \quad \text{with} \quad \| Q \| \leq 2,  
\end{displaymath} 
we have $b_n (D) \leq  b_n (M) \, \| Q \| \leq 2 \, b_n (M)$, that is:
\begin{displaymath} 
a_n (S_\mu) = b_n (S_\mu) \geq b_n (M) \geq \frac{b_n (D)}{2} = \frac{m_{n, n}}{2} \geq \frac{\varepsilon_n^{2}}{64} \, \raise 1 pt \hbox{,}
\end{displaymath} 
since the $n$ first approximation numbers of the diagonal matrix $D$ (the matrices being viewed as well as operators on the Hilbertian space $\C^n$ with its 
canonical basis) are $m_{1, 1}, \ldots, m_{n, n}$. It follows that, using \eqref{equal}:
\begin{equation}\label{injection} 
a_n (I_\mu) = a_n (I_\mu^{*}) = \sqrt{a_n (S_\mu)} \geq \frac{\varepsilon_{n}}{8} \, \cdot
\end{equation}
In view of \eqref{equal}, we have as well $a_n (C_\varphi) \geq \eps_n / 8$, and we are done. \qed 
\par\medskip

\noindent {\bf Proof of Lemma~\ref{cle}.} Write  $M = (m_{i, j}) = D (I + N)$ with $N = D^{- 1}R$. One has:
\begin{equation}\label{aussi} 
N = (\nu_{i, j}), \quad \text{with} \quad \nu_{i, i} = 0 \quad \text{and} \quad \nu_{i, j} = \frac{m_{i, j}}{m_{i,i}} \text{ for } j \neq i \, .
\end{equation}
We shall show that $\Vert N \Vert \leq 1/2$ by using the (unweighted) Schur test, which we recall (\cite{Halmos-livre}, Problem~45):
\begin{proposition}\label{test} 
Let $(a_{i, j})_{1 \leq i, j \leq n}$ be a matrix of complex numbers. Suppose that there exist two positive numbers $\alpha, \beta > 0$ such that: \par \smallskip 

$1.$ $\sum_{j = 1}^n \vert a_{i, j} \vert \leq \alpha $ for all $i$; \par \smallskip 

$2.$ $\sum_{i = 1}^n \vert a_{i, j} \vert \leq \beta $ for all $j$. \par\smallskip

\noindent Then, the (Hilbertian) norm of this matrix  satisfies $\Vert A \Vert \leq \sqrt{\alpha \beta}$. 
\end{proposition}
\par\smallskip 

It is essential for our purpose to note  that: 
\begin{align} 
& i < j \quad \Longrightarrow  \quad \vert \nu_{i, j} \vert \leq 32 \, \delta^{j - i} \, , \label{useful} \\
& i > j \quad \Longrightarrow \quad \vert \nu_{i, j} \vert \leq 32 \, (2 \, \delta)^{i - j} \label{useful2} \, .
\end{align}
Indeed, we see from \eqref{dimanche} and \eqref{aussi} that, for $i < j$:
\begin{displaymath} 
\vert \nu_{i, j}\vert = \frac{\vert m_{i, j} \vert}{m_{i, i}} \leq 32\, \eps_i \, \eps_j \, \eps_i^{- 2} \delta^{j - i} \leq 32 \, \delta^{j - i}  
\end{displaymath} 
since  $\eps_j \leq \eps_i$. Secondly, using $\varepsilon_j/ \varepsilon_i \leq 2^{i - j}$ for $i > j$ (recall that we assumed that 
$\varepsilon_{k + 1} \geq \varepsilon_k / 2$), as well as $\vert m_{i, j}\vert = \vert m_{j, i} \vert$, we have, for $i > j$: 
\begin{displaymath} 
\vert \nu_{i, j} \vert = \frac{\vert m_{j, i} \vert}{m_{i, i}} \leq 32 \, \frac{\eps_j}{\eps_i} \, \delta^{i - j} \leq 32 \, (2\, \delta)^{i - j} .
\end{displaymath} 

Now, for fixed $i$, \eqref{useful} gives:
\begin{align*}
\sum_{j = 1}^n \vert \nu_{i, j}\vert 
& = \sum_{j > i} \vert \nu_{i, j} \vert + \sum_{j < i} \vert \nu_{i, j} \vert 
\leq 32 \,\bigg( \sum_{j > i} \delta^{j - i} + \sum_{j < i} (2\, \delta)^{i - j} \bigg) \\ 
& \leq 32 \bigg( \frac{\delta}{1 - \delta} + \frac{2\, \delta}{1 - 2 \, \delta} \bigg) \leq 32 \, \frac{3\, \delta}{1 - 2 \, \delta} 
\leq \frac{96}{198} \leq \frac{1}{2} \, ,
\end{align*}
since $\delta \leq 1/200$. Hence:
\begin{equation}\label{first} 
\sup_{i}\Big( \sum_{j} \vert \nu_{i, j} \vert \Big) \leq 1/2 \, .
\end{equation}
In the same manner, but using \eqref{useful2} instead of \eqref{useful}, one has:
\begin{equation} \label{second} 
\sup_{j} \Big( \sum_{i} \vert \nu_{i, j} \vert \Big) \leq 1/2 \, . 
\end{equation}
Now, \eqref{first}, \eqref{second} and the  Schur criterion recalled above give:
\begin{displaymath} 
\Vert N \Vert\leq \sqrt{1/2 \times 1/2} = 1/2 \, , 
\end{displaymath} 
as claimed. \qed 
\bigskip

\noindent {\bf Remark.} We could reverse the point of view in the preceding proof: start from $\theta$ and see what lower bound for $a_{n} (C_\varphi)$ 
emerges. For example, if $\theta (h) \approx h$ as is the case for lens maps (see \cite{PDHL}), we find again that 
$a_n (C_\varphi) \geq \delta_0 > 0$ and that $C_\varphi$ is not compact. But if $\theta (h) \approx h^{1 + \alpha}$ with $\alpha > 0$, the method only gives 
$a_n (C_\varphi) \gtrsim \e^{- \alpha n}$ (which is always true: see \cite{PDHL}, Theorem~2.1), whereas the methods of \cite{PDHL} easily give 
$a_n (C_\varphi) \gtrsim \e^{- \alpha \sqrt{n}}$. Therefore, this $\mu$-method seems to be sharp when we are close to non-compactness, and to be beaten by 
those of \cite{PDHL} for ``strongly compact'' composition operators. 


\subsection{Optimality of the EKSY result} \label{proof Th 2} 

El Fallah, Kellay, Shabankhah and Youssfi proved in \cite{EKSY} the following: if $\varphi$ is a Schur function such that $\varphi \in \mathcal{D}$ and 
$\Vert \varphi^p \Vert_{\mathcal{D}}  = O\, (1)$ as $p \to \infty$, then $\varphi$ is a symbol on $\mathcal{D}$. 
We have the following theorem, already stated  in the Introduction, which shows the optimality of their result.

\begin{theorem} \label{mieux} 
Let $(M_p)_{p \geq 1}$ be an arbitrary sequence of positive numbers such that $\lim_{p \to \infty} M_p = \infty$. Then, there exists a Schur function 
$\varphi \in \mathcal{D}$ such that: \par\smallskip

$1)$ $\Vert \varphi^p \Vert_{\mathcal{D}} = O\, (M_p)$ as $p \to \infty$; \par\smallskip 

$2)$ $\varphi$ is not a symbol on $\mathcal{D}$.
\end{theorem}

\noindent{\bf Remark.} We first observe that we cannot replace $\lim$ by $\limsup$ in Theorem~\ref{mieux}. Indeed, since $\varphi \in \mathcal{D}$, the 
measure $\mu = n_{\varphi} \, dA$ is finite, and
\begin{displaymath} 
\| \varphi^p \|_{\mathcal{D}}^2 = p^2 \int_{\D} | w |^{2 p - 2} \, d\mu (w) 
\geq c \, p^2 \bigg( \int_{\D} | w |^2 \, d\mu (w) \bigg)^{p - 1} \geq c \, \delta^p \, , 
\end{displaymath} 
where $c$ and $\delta$ are positive constants.\par\medskip 

\noindent {\bf Proof of Theorem~\ref{mieux}.} We may, and do, assume that $(M_p)$ is non-decreasing and integer-valued. Let $(l_n)_{n \geq 1}$ be an 
non-decreasing sequence of positive integers tending to infinity, to be adjusted. Let $\Omega$ be  the subdomain of the right half-plane $\C_0$  defined as follows. 
We set: 
\begin{displaymath} 
\varepsilon_n = - \log (1 - 2^{-n}) \sim 2^{- n} \, , 
\end{displaymath} 
and we consider the (essentially) disjoint boxes ($k = 0, 1, \ldots$): 
\begin{displaymath} 
B_{k, n} = B_{0, n} + 2 k \pi i \, ,
\end{displaymath} 
with:
\begin{displaymath} 
B_{0, n} = \{u \in \C \, ; \ \varepsilon_{n + 1} \leq \Re u \leq \varepsilon_n \text{ and } | \Im u | \leq 2^{- n} \pi \} \, ,
\end{displaymath} 
as well as  the union 
\begin{displaymath} 
T_n = \bigcup_{0 < k < l_n} B_{k, 2n} \, , 
\end{displaymath} 
which is a kind of broken tower above the "basis" $B_{0, 2n}$ of even index. \par 

We also consider, for $1 \leq k \leq l_n - 1$, very thin vertical pipes $P_{k, n}$ connecting $B_{k, 2n}$ and $B_{k - 1, 2n}$, of side lengths $4^{- 2n}$ and 
$2 \pi (1 - 2^{- 2n})$ respectively:
\begin{displaymath} 
P_{k, n} = P_{0, n} + 2 k \pi i \, ,
\end{displaymath} 
and we set:
\begin{displaymath} 
P_n = \bigcup_{1 \leq k < l_n} P_{k, n}
\end{displaymath} 

Finally, we set: 
\begin{displaymath} 
F = \left(\bigcup_{n = 2}^\infty B_{0, n} \right) \cup \left( \bigcup_{n = 1}^\infty T_n \right)  \cup \left( \bigcup_{n = 1}^\infty P_n \right)
\end{displaymath} 
and:
\begin{displaymath} 
\Omega = \mathop{F}^{\circ}
\end{displaymath} 
\begin{figure} [h]
\centering
\includegraphics[width=5cm]{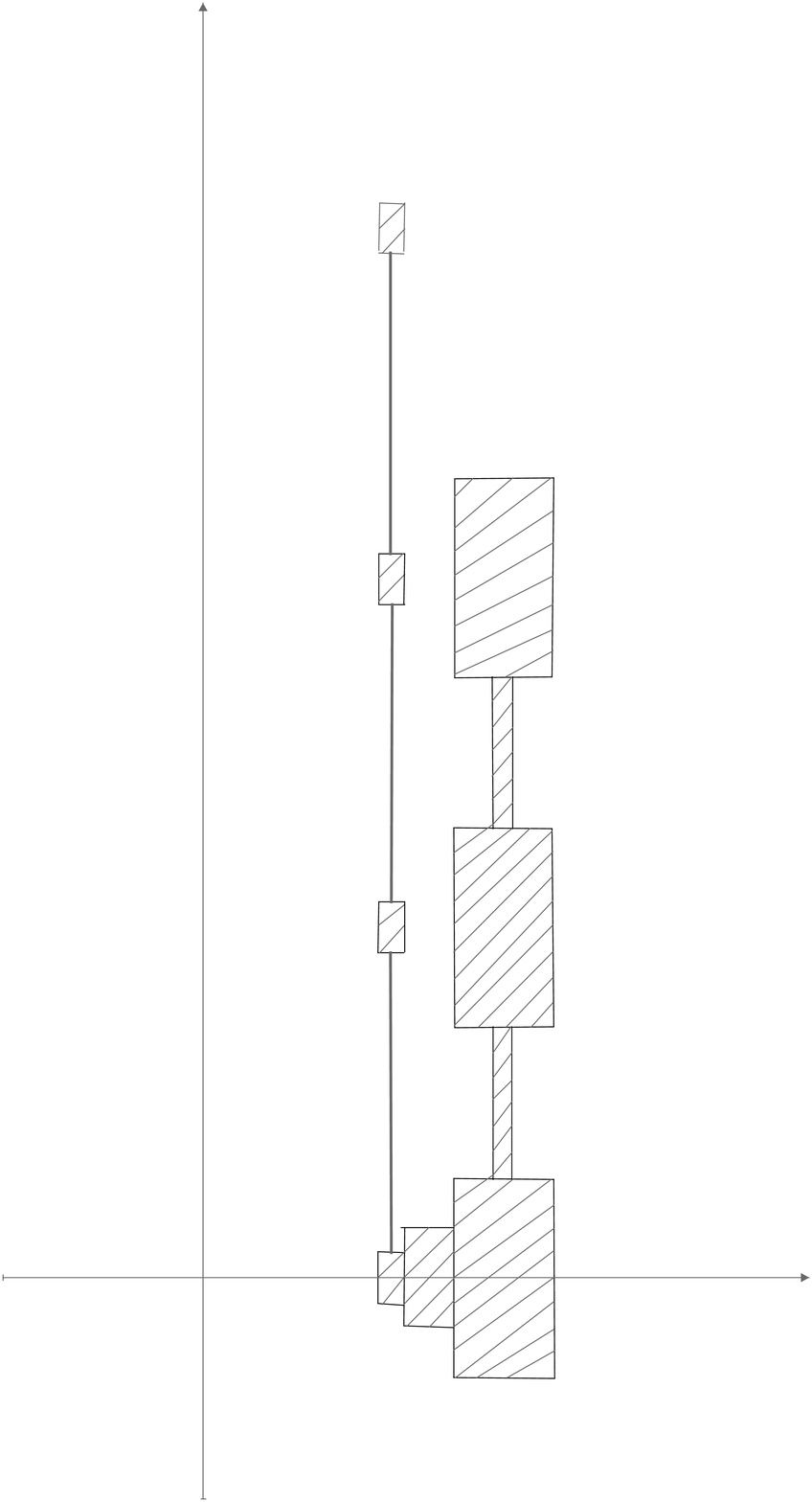}
\end{figure}

Then $\Omega$ is a simply connected domain. Indeed, it is connected thanks to the $B_{0, n}$ and the $P_n$, since the $P_{k, n}$ were added to ensure that. 
Secondly, its unbounded complement is connected as well, since we take one value of $n$ out of two in the union of sets $B_{k, n}$ defining $F$. \par\smallskip 
 
Let now $f \colon \D \to \Omega$ be a Riemann map, and $\varphi = \e^{- f} \colon \D \to \D$. \par\medskip 

We introduce the Carleson window $W = W (1, h)$ defined as: 
\begin{displaymath} 
W (1, h) = \{z \in \D \, ; \  1 - h \leq |z| < 1 \text{ and } | \arg z | < \pi\,h \} \, .
\end{displaymath} 
This is a variant of the sets $S (1, h)$ of Section~\ref{notations}. We also introduce  the Hastings-Luecking half-windows $W'_n$ defined by:
\begin{displaymath} 
W'_n = \{z \in \D \, ; \ 1 - 2^{- n} < | z | < 1 - 2^{- n - 1} \text{ and } | \arg z | < \pi \, 2^{- n}\}.
\end{displaymath} 
We will also need the sets:
\begin{displaymath} 
E_n = \e^{- (T_n \cup B_{0, 2 n + 1} \cup P_n)} = \e^{- (B_{0, 2n} \cup B_{0, 2 n + 1} \cup P_{0, n})} \, ,
\end{displaymath} 
for which one has:
\begin{displaymath} 
\varphi (\D) \subseteq \bigcup_{n = 1}^\infty E_n \, .
\end{displaymath} 

Next,  we  consider the measure $\mu = n_{\varphi} \, dA$, and a Carleson window $W = W (1, h)$ with $h = 2^{- 2N}$. We observe that 
$W'_{2N} \subseteq W$ and claim that:

\begin{lemma} \label{last lemma}
One has: \par \smallskip

$1)$ $w \in W'_{2N} \quad \Longrightarrow \quad n_\phi (w) \geq l_N$; \par\smallskip 

$2)$ $\| \varphi^p \|_{\mathcal{D}}^2 \lesssim p^2 \sum_{n = 1}^\infty l_n \, 16^{- n} \, \e^{- p \, 4^{- n}}$.
\end{lemma}
\noindent{\bf Proof of Lemma~\ref{last lemma}.} $1)$ Let $w = r\, \e^{i \theta} \in W'_{2N}$ with $1 - 2^{- 2N} < r < 1 - 2^{- 2 N - 1}$ and 
$|\theta| < \pi \, 2^{- 2N}$. As $- (\log r + i \theta) \in B_{0, 2N}$, one has $- (\log r + i \theta) = f (z_0)$ for some $z_0 \in \D$. Similarly, 
$- (\log r + i \theta) + 2 k \pi i$, for $1 \leq k < l_N$, belongs to $B_{k, 2N}$ and can be written as $f (z_k)$, with $z_k \in \D$. The $z_k$'s, $0 \leq k < l_N$, 
are distinct and satisfy $\varphi (z_k) = \e^{- f (z_k)} = \e^{- f (z_0)} = w$ for $0 \leq k < l_N$, thanks to the $2\pi i$-periodicity of the exponential function. 
\par\smallskip

$2)$ We  have $A (E_n) \lesssim \e^{- 2 \varepsilon_{2 n + 2}} 4^{- 2 n} \leq 4^{- 2 n}$ (the term $\e^{- 2 \varepsilon_{2 n + 2}}$ coming from the Jacobian 
of $\e^{- z}$) and we  observe that 
\begin{displaymath} 
w \in E_n \quad \Longrightarrow \quad | w |^{2 p - 2} \leq (1 - 2^{- 2 n - 1})^{2 p - 2} \lesssim \e^{- p\, 4^{- n}} \, .
\end{displaymath}   
It is easy to see that $n_\phi (w) \leq l_n$ for $w \in E_n$; thus we obtain, forgetting the constant term $|\phi (0)|^{2p} \leq 1$, using \eqref{utile} and keeping 
in mind the fact that $n_\phi (w) = 0$ for $w \notin \phi (\D)$: 
\begin{align*}
\Vert \varphi^p\Vert_{\mathcal{D}}^{2} 
& = p^2 \int_{\varphi(\D)} \vert w \vert^{2 p - 2} \, n_\varphi (w) \, dA (w) \\ 
& \leq p^2 \bigg( \sum_{n = 1}^\infty \int_{E_n} \vert w \vert^{2 p - 2} \, n_\varphi (w) \, dA (w) \bigg) \\
& \leq p^2 \bigg( \sum_{n = 1}^\infty  \int_{E_n} \vert w \vert^{2 p - 2} \, l_n \, dA (w) \bigg) \\
& \lesssim p^2 \sum_{n = 1}^\infty l_{n} \, 16^{- n} \, \e^{- p \, 4^{- n}} \, ,
\end{align*}
ending the proof of Lemma~\ref{last lemma}. \qed 
\par\bigskip

\noindent \emph{End of the proof of Theorem~\ref{mieux}.} Note that, as a consequence of the first part of the proof of Lemma~\ref{last lemma}, one has 
\begin{displaymath} 
\mu (W) \geq \mu (W'_{2 N}) = \int_{W'_{2 N}} n_{\varphi} \, dA \geq l_N A (W'_{2 N}) \gtrsim l_N h^2 \, ,
\end{displaymath} 
which implies that $\sup_{0 < h < 1} h^{- 2} \mu [W (1, h)] = + \infty$ and shows that $C_\varphi$ is not bounded on $\mathcal{D}$ by Zorboska's criterion 
(\cite{ZOR}, Theorem~1), recalled in \eqref{ninaz}.  \par

It remains now to show that we can adjust the non-decreasing sequence of integers $(l_n)$ so as to have $\Vert \varphi^p \Vert_{\mathcal{D}} = O\, (M_p)$. 
To this effect, we first observe that, if one sets $F (x) = x^2 \, \e^{- x}$, we have:
\begin{displaymath} 
p^2 \sum_{n = 1}^\infty  16^{- n} \, \e^{- p \, 4^{- n}} = \sum_{n = 1}^\infty F \left(\frac{p}{4^n} \right) \lesssim 1 \, .
\end{displaymath} 
Indeed, let $s$ be the integer such that $4^s \leq p < 4^{s+1}$. We have:
\begin{displaymath}
\sum_{n = 1}^\infty  F \left(\frac{p}{4^n} \right) 
\lesssim \sum_{n = 1}^s  \frac{4^{n}}{p} + \sum_{n > s}  F (4^{- (n - s - 1)})
\lesssim 1 + \sum_{n = 0}^\infty F (4^{- n}) < \infty \, ,
\end{displaymath}
where we used that $F$ is increasing on $(0, 1)$ and satisfies $F (x) \lesssim \min (x^2, 1/x)$ for $x > 0$. We finally choose the non-decreasing sequence 
$(l_n)$ of integers as:
\begin{displaymath} 
l_n = \min (n, M_n^{2} ) \, .
\end{displaymath} 
In view of Lemma~\ref{last lemma} and of the previous observation, we obtain: 
\begin{align*}
\Vert \varphi^p \Vert_{\mathcal{D}}^2 
& \lesssim p^2 \sum_{n = 1}^\infty  16^{- n} \, \e^{- p \, 4^{- n}} l_n  \\
& \leq p^2 \sum_{n = 1}^p 16^{- n} \, \e^{- p \, 4^{- n}} l_p + p^2 \sum_{n > p} 16^{- n} l_n \\ 
& \lesssim l_p + p^2 \sum_{n > p} 4^{- n} 
\lesssim l_p + p^2 \, 4^{- p} \lesssim M_p^2 \, ,
\end{align*}
as desired. This choice of $(l_n)$  gives us an unbounded composition operator on $\mathcal{D}$ such that $\| \phi^p \|_{\mathcal{D}} = O\, (M_p)$, which 
ends the proof of Theorem~\ref{mieux}. \qed


\bigskip

\noindent
{\rm Daniel Li}, Univ Lille Nord de France, \\
U-Artois, Laboratoire de Math\'ematiques de Lens EA~2462 \\ 
\& F\'ed\'eration CNRS Nord-Pas-de-Calais FR~2956, \\
Facult\'e des Sciences Jean Perrin, Rue Jean Souvraz, S.P.\kern 1mm 18, \\
F-62\kern 1mm 300 LENS, FRANCE \\ 
daniel.li@euler.univ-artois.fr
\medskip

\noindent
{\rm Herv\'e Queff\'elec}, Univ Lille Nord de France, \\
USTL, Laboratoire Paul Painlev\'e U.M.R. CNRS 8524 \& 
F\'ed\'eration CNRS Nord-Pas-de-Calais FR~2956, \\
F-59\kern 1mm 655 VILLENEUVE D'ASCQ Cedex, 
FRANCE \\ 
Herve.Queffelec@univ-lille1.fr
\smallskip

\noindent
{\rm Luis Rodr{\'\i}guez-Piazza}, Universidad de Sevilla, \\
Facultad de Matem\'aticas, Departamento de An\'alisis Matem\'atico \& IMUS,\\ 
Apartado de Correos 1160,\\
41\kern 1mm 080 SEVILLA, SPAIN \\ 
piazza@us.es\par

\end{document}